\documentclass[10pt]{article}
	\usepackage{amsmath}
      \usepackage{amssymb}
      \usepackage{pxfonts}
      \usepackage{physics}
        \usepackage{enumitem}
        \usepackage{ragged2e}
        \usepackage{graphicx}
        \usepackage{caption}
        \usepackage{subcaption}
        \usepackage{float}
        \usepackage{mathrsfs}
        \usepackage{mathtools}
        \usepackage{hyperref}
        \usepackage[dvipsnames]{xcolor}
        \usepackage{cite} 
        \usepackage[margin=1in, top=1in, bottom=1in]{geometry}

        \numberwithin{equation}{section}
        \newtheorem{theorem}{Theorem}
        \usepackage{authblk}
	
\begin{document}
	
        \justifying
        \title{Global Dynamics and Patterns of a Gray--Scott Model with Local--Nonlocal Diffusion}
        \author{Md Shah Alam\thanks{Email: \href{mailto:malam@htu.edu}{malam@htu.edu}}}
\affil{Department of Mathematics, Huston--Tillotson University, Austin, TX 78702, USA.}
%\affil{Department of Mathematics, University of Houston, Houston, TX 77004, USA.}

\date{}
\maketitle

\begin{abstract}

\noindent
This paper studies the global existence of component-wise nonnegative solutions of the Gray-Scott model in $\Omega \subset \mathbb{R}^n$, $n \ge 1$, with a mixture of both local and nonlocal diffusion operators. We use semigroup theory with duality arguments to establish the global existence and boundedness of solutions of our model. We also present the patterns of both local and mixed Gray-Scott models.
\end{abstract}



\section{Introduction and Main Results}\label{Introduction and Main Results}
\subsection{Introduction}
The classical Gray--Scott model equipped with local (Laplacian) diffusion has been extensively investigated due to its ability to generate a wide repertoire of self-organized spatial and spatiotemporal patterns, as well as its analytical tractability. A substantial body of work has established the existence, boundedness, stability, and bifurcation structure of solutions across spatial dimensions and parameter regimes; see, for instance, \cite{nishiura2001spatio,doelman1997pattern,mcgough2004pattern,muratov2001spike,doelman1998stability,ali2023spatiotemporal}. The classical local Gray--Scott system takes the form
\begin{flalign}\label{eq:local}
\begin{cases}
u_t(x,t) = d_1 \Delta u(x,t) - u(x,t) v(x,t)^2 + f(1- u(x,t)), & (x,t) \in \Omega \times (0,\infty),\\[0.3em]
v_t(x,t) = d_2 \Delta v(x,t) + u(x,t) v(x,t)^2 - (f+ \kappa)v(x,t), & (x,t) \in \Omega \times (0,\infty),\\[0.3em]
\frac{\partial u}{\partial \eta} = \frac{\partial v}{\partial \eta} = 0, & (x,t) \in \partial \Omega \times (0,\infty),\\[0.3em]
u(x,0)=u_0(x)\ge0,\quad v(x,0)=v_0(x)\ge0, & x \in \Omega,
\end{cases}
\end{flalign}
where $u$ and $v$ denote concentrations of two reacting species, $d_1,d_2> 0$ are diffusion coefficients, and $f, \kappa >0$ represent the feed and removal rates, respectively. The operator $\frac{\partial}{\partial \eta}$ denotes the outward normal derivative along $\partial \Omega$, corresponding to homogeneous Neumann boundary conditions.

\vspace{0.2em}
\noindent
The Gray--Scott model provides a minimal yet remarkably expressive prototype for studying self-organization in reaction--diffusion systems. Unlike classical activator--inhibitor mechanisms---in which one species promotes growth while the other suppresses it---the Gray--Scott dynamics are driven by an autocatalytic reaction in which the activator both amplifies itself and depletes the inhibitor. This feedback structure is responsible for the emergence of spots, stripes, labyrinths, spiral waves, and spatiotemporal chaos observed in experiments and numerical simulations.

\vspace{0.2em}
\noindent
From the analytical perspective, the foundational works of Hollis, Martin, and Pierre~\cite{hollis1987global} established global existence and boundedness for two-species reaction--diffusion systems under structural growth conditions on the nonlinearities. Morgan~\cite{morgan1989global} subsequently generalized these results to systems with an arbitrary number of components $m \ge 2$, providing a robust framework for global analysis of multi-species diffusion systems.

\vspace{0.2em}
\noindent
More recently, motivated by applications in ecology, biology, and material science, attention has shifted to \emph{nonlocal} diffusion operators that model long-range dispersal or interaction effects. For variants of the Gray--Scott model in which one of the diffusion operators is replaced by a convolution-type nonlocal operator, global existence results have been obtained by Laurencot and Walker~\cite{laurencot2023nonlocal} for general spatial dimensions $n \ge 1$, and by Cappanera, Jaramillo, and Ward~\cite{cappanera2024analysis} in the one-dimensional case. These works demonstrate that nonlocality fundamentally alters the spatial dynamics and stability properties of the system.
\noindent
Nonlocal diffusion naturally arises in contexts where transport occurs across multiple spatial scales, such as organism dispersal in ecology, chemical exchange in layered media, and cell-to-cell signaling in developmental biology. In contrast, local diffusion corresponds to short-range Brownian motion. Models incorporating both mechanisms capture richer physical behavior and are increasingly relevant in multiscale biological and material systems.

\vspace{0.2em}
\noindent
In this work, we study the system \eqref{eq:1.1}, where the activator $u$ undergoes nonlocal diffusion mediated by a convolution operator, while the inhibitor $v$ diffuses locally through the Laplacian. This configuration is of particular interest: the nonlocal diffusion of $u$ models long-range aggregation or transport, whereas the local diffusion of $v$ captures short-range spread. We emphasize that the placement of the nonlocal operator is crucial. If the roles of the diffusion operators are reversed---that is, if $u$ is governed by the Laplacian and $v$ by the nonlocal operator---then the solutions may develop finite-time blow-up, even in low spatial dimensions.

\vspace{0.2em}
\noindent
The novelty of our work lies in the hybrid diffusion structure of the system. While both local and nonlocal Gray--Scott models have been studied independently, the coupled interaction of these \emph{two distinct diffusion mechanisms} within a two-dimensional reaction--diffusion framework has received very limited analytical attention. The presence of nonlocal diffusion in only one equation breaks the symmetry of the system and prevents the direct application of classical energy methods, maximum principles, and comparison arguments typically used in the study of reaction--diffusion systems. Consequently, new analytical techniques are required to obtain global existence and uniform bounds, especially in two dimensions where nonlinear effects are strongest.

\vspace{0.2em}
\noindent
Our main contribution is to establish global existence, uniform boundedness, and regularity of classical solutions to \eqref{eq:1.1} under minimal structural assumptions on the kernel and nonlinear reaction terms. Our analysis employs the semigroup approach for parabolic systems, together with a duality method that compensates for the lack of direct comparison principles in the presence of nonlocal diffusion. This strategy yields uniform a priori bounds and leads to the global existence and boundedness of classical solutions to \eqref{eq:1.1}. Our results extend the classical theory for local diffusion systems to a novel hybrid local–nonlocal Gray--Scott framework and provide analytical support for complex pattern formation observed in numerical simulations.

\subsection{Problem Settings:}
\noindent
Let $n \in \mathbb{N}$ and $\Omega$ be a bounded open subset of $\mathbb{R}^n$ with smooth boundary $\partial \Omega$ ($\partial \Omega$ is an $n-1$ dimensional $C^{2+\mu}$ $(\mu \in (0,1))$ manifold such that $\Omega$ lies locally on one side of $\partial \Omega$).\\
\noindent
Our goal is to study the following nonlocal Gray-Scott model equipped with both local and nonlocal convolution diffusion operator which shows anomalous faster diffusion
\begin{flalign}\label{eq:1.1}
\begin{cases}
u_t(x,t) & = d_1 \Gamma u(x,t) - u(x,t) v(x,t)^2 + f(1- u(x,t)), \hspace{0.37 in} (x,t) \in \Omega \times (0,\infty),\\
v_t(x,t) & = d_2 \Delta v(x,t) + u(x,t) v(x,t)^2 - (f+ \kappa)v(x,t), \hspace{0.36 in} (x,t) \in \Omega \times (0,\infty),\\
\frac{\partial v}{\partial \eta} & = 0, \hspace{2.84 in} (x,t) \in \partial \Omega \times (0,\infty),\\
u(x,0) & = u_0(x) \ge 0, \ v(x,0) =v_0(x) \ge 0\hspace{1.16  in} x \in \Omega.
\end{cases}
\end{flalign}
\noindent
We consider nonnegative initial data \(u_0 \in C(\overline{\Omega})\) and \(v_0 \in C^2(\overline{\Omega})\). The diffusion coefficients are distinct and positive, and the parameters \(f,\kappa>0\) retain the same meaning as in \eqref{eq:local}. System \eqref{eq:1.1} is supplemented with homogeneous Neumann boundary conditions, as in the classical local Gray--Scott model. The nonlocal diffusion operator \(\Gamma\) acting on a function \(z:\Omega \to \mathbb{R}\) is defined by
\begin{equation}\label{eq:1.2}
\Gamma z(x) := \int_{\Omega} \phi(x,y)\,[z(y)-z(x)]\,dy, \qquad x \in \Omega,
\end{equation}
where the kernel \(\phi : \overline{\Omega} \times \overline{\Omega} \to [0,\infty)\) is measurable, nonnegative, continuous, and symmetric, i.e.,
\[
\phi(x,y) = \phi(y,x), \qquad \text{for all } x,y \in \Omega,
\]
and satisfies the mass bound
\begin{equation}\label{eq:1.3}
\int_{\Omega} \phi(x,y)\,dy \le \lambda < \infty, \qquad \text{for all } x \in \Omega.
\end{equation}
\noindent
We define the reaction terms
\begin{equation*}
g_1(u,v) = -uv^2 + f(1-u), \qquad g_2(u,v) = uv^2 - (f+\kappa)v.
\end{equation*}
Let us assume that
\begin{enumerate}
  \item[(D)] \hypertarget{item:D} (Diffusive constants) $d_1,d_2> 0$;
  \item[(f1)] \hypertarget{item:f1} (Local Lipschitz continuity) $g_1,g_2 : \mathbb{R}^2  \times \mathbb{R}_+  \to \mathbb{R}^2$ is locally Lipschitz with respect to $u,v$ and uniformly Lipschitz with respect to $t$;
  \item[(f2)] \hypertarget{item:f2} (Quasi-positivity) For each $g_1(u,0),g_2(0,v) \ge 0$ whenever $t \ge 0$ and $u,v \in \mathbb{R}_+^2$ such that $u,v=0$.
\end{enumerate}
\hyperlink{item:f2}{(f2)} ensures the componentwise nonnegativity of the solution along with nonnegative initial data.

\subsection{Statement of the Main Result}
\begin{theorem}\label{Th:1}
Assume that (\hyperlink{item:D}{(D)}, \hyperlink{item:f1}{(f1)}, \hyperlink{item:f2}{(f2)}, \eqref{eq:1.2} and \eqref{eq:1.3} hold, and let 
$u_0 \in C(\overline{\Omega}, \mathbb{R}_+)$ and 
$v_0 \in C^2(\overline{\Omega}, \mathbb{R}_+)$ satisfy 
$\dfrac{\partial v_0}{\partial \eta} = 0$ on $\partial \Omega$. 
Then there exists a unique classical global componentwise nonnegative solution $(u(x,t), v(x,t))$ to \eqref{eq:1.1}, and both $u$ and $v$ remain uniformly bounded in the sup-norm for all time.
\end{theorem}
\noindent
\textbf{Organization of the paper.}
Section~\ref{Proof of Main Result} is devoted to the proof of our main result. Section~\ref{sec:patterns} presents numerical simulations comparing pattern formation in the classical and hybrid Gray--Scott models. In Section~\ref{Conclusion}, we summarize our findings and discuss potential avenues for future research.

\section{Proof of Main Result}\label{Proof of Main Result}
Define $X := C(\overline\Omega)$ with $X_+ := \{ z \in X: z \ge 0 \ \text{in } \Omega\}$. Then using (\ref{eq:1.2}) and (\ref{eq:1.3}), for any $x \in \Omega$ we can write,   
\begin{flalign*}
 \vert \Gamma_{\phi}z(x) \vert & \le \int_{\Omega} \phi(x,y) \vert z(y)-z(x) \vert \ dy \le \int_{\Omega} \phi(x,y) \vert z(y) \vert \ dy +\int_{\Omega} \phi(x,y) \vert -z(x) \vert \ dy \le 2 \lambda \Vert z  \Vert_{\infty}.
\end{flalign*}
\noindent
As a result, $T_1(t)=e^{\Gamma t}$ is a strongly continuous semigroup on $X$. We also define $A(z)= -d_2 \Delta z$ with the following domain,
\begin{flalign*}
z \in D(A) & = \big\{w\in W_p^2(\Omega)\,: p>n, \ A(w)\in C(\overline\Omega),\,\frac{\partial w}{\partial\eta}=0\,\text{on }\partial\Omega \big\}.
\end{flalign*}
Then $A:D(A)\to X$ is an unbounded linear operator that generates an analytic semigroup on $X$ \cite{stewart1974generation}, and hence a strongly continuous semigroup $T_2(t)$ on $X$.\\
\\
\noindent
\textbf{Proof of Theorem {\ref{Th:1}}:} With the above settings, results in \cite{pazy2012semigroups} guarantee that (\ref{eq:1.1}) has a unique maximal solution $(u,v)$ on a maximal interval $[0,T_m)$, where $T_m$ is a positive extended real number, given by 
\begin{flalign*}
    u(x,t) & = T_1(t) u_0(x) +\int_0^t T_1(t-\tau) g_1(u(x,\tau),v(x,\tau))d\tau \\
    v(x,t) & = T_2(t) v_0(x) +\int_0^t T_2(t-\tau) g_2(u(x,\tau),v(x,\tau))d\tau
\end{flalign*}
It is straightforward to prove that \textit{quasi-positivity} implies the maximal solution $(u,v)$ is component-wise nonnegative. It is well known that if $T_m<\infty$ then
\begin{flalign}
  \limsup_{t \to T_m^-} \big\{ \Vert u(\cdot,t)\Vert_{\infty,\Omega} + \Vert v(\cdot,t)\Vert_{\infty,\Omega} \big\} = \infty. \label{eq:2.1}
\end{flalign}
Consequently, we demonstrate global existence by proving solutions cannot blow up in the sup norm in finite time.
Now, if we define $A_1(u)=d_1 \Gamma u-fu$, then $A_1(u)$ generates $T_{11}(t)$, and from results in \cite{laurencot2023nonlocal} we observe that $\Vert T_{11}(z) \Vert_{\infty,\Omega} \le e^{-tf}\Vert z \Vert_{\infty,\Omega}$ for all $z \in C(\overline\Omega,\mathbb{R})$. In addition,
\begin{flalign}
u_t(x,t) & = d_1 \Gamma u(x,t) - u(x,t) v(x,t)^2 + f(1- u(x,t)) \nonumber\\
=> \ u(x,t) & \le T_{11}(t) u_0(x) + \int_0^t T_{11}(t-\tau) f \ d\tau \hspace{1.4 in} [\because \ u,v\ge 0 \ \forall \ t>0] \nonumber\\
=> \ u(x,t) & \le \Vert T_{11}(t) u_0 \Vert_{\infty,\Omega} + \int_0^t \Vert T_{11}(t-\tau) f \Vert_{\infty,\Omega} \ d\tau \le e^{-tf} \Vert u_0 \Vert_{\infty,\Omega} + 1 - e^{-tf} \nonumber\\
\therefore \ \Vert u(\cdot,t) \Vert_{\infty,\Omega} & \le C , \ \ \forall t>0. \hspace{2 in} \left[ C:=\max\big\{\Vert u_0 \Vert_{\infty,\Omega}, 1 \big\}\right]  \label{eq:2.2}
\end{flalign}
\noindent
Now, adding the first two equations of the (\ref{eq:1.1}) we get,
\begin{flalign*}
  u_t(x,t) + v_t(x,t) & = d_1 \Gamma u(x,t) + d_2 \Delta v(x,t) + f(1- u(x,t)) -(f+\kappa) v(x,t)\\
=> \ \frac{\partial}{\partial t} \big(u(x,t) + v(x,t)\big) & \le d_1 \Gamma u(x,t) + d_2 \Delta v(x,t) + f - f\tilde{\kappa}(u(x,t)+v(x,t)) \hspace{0.7 cm} [\tilde{\kappa}=\min\big\{ \kappa,1 \big\}>0]\\
=> \ \frac{d}{dt} \int_{\Omega} \big(u(x,t) + v(x,t)\big) dx & \le d_1 \int_{\Omega} \Gamma u(x,t) dx + d_2 \int_{\Omega} \Delta v(x,t) dx + f \int_{\Omega} dx - f\tilde{\kappa} \int_{\Omega} (u(x,t)+v(x,t)) dx
\end{flalign*}
where the boundary condition gives $\int_{\Omega} \Delta v(x,t) dx=0$. Note that $\phi(x,y)=\phi(y,x)$ immediately implies that
\begin{flalign*}
d_1 \int_{\Omega} \Gamma u(x,t) dx = d_1 \int_{\Omega} \int_{\Omega} \phi(x,y) (u(y,t)-u(x,t))dxdy = 0.
\end{flalign*}
Therefore, we get,
\begin{flalign}
\frac{d}{dt} \int_{\Omega} \big(u(x,t) + v(x,t)\big) dx & \le f\vert \Omega \vert - f\tilde{\kappa} \int_{\Omega} (u(x,t)+v(x,t)) dx \nonumber\\
\therefore \ \int_{\Omega} \big(u(x,t) + v(x,t)\big) dx & \le M\ \ \ \forall t>0. \ \ \ \ \ [M = \max\big\{\frac{\vert \Omega \vert}{\tilde{\kappa}},\int_{\Omega}(u_0(x)+v_0(x))dx \big\}]\label{eq:2.3}
\end{flalign}
By (\ref{eq:2.2}) and the boundedness of $\Omega$, we have, for every $1\le p<\infty$ and $0\le \tau<T<T_m$,
\begin{flalign}\label{eq:u_p}
\|u\|_{p,\Omega\times(\tau,T)} & =
\left(\int_{\tau}^{T}\int_{\Omega}|u(x,t)|^p\,dxdt\right)^{\frac1p} \le |\Omega|^{\frac1p}(T-\tau)^{\frac1p} \|u\|_{\infty,\Omega\times(\tau,T)} \le M_p(T-\tau),
\end{flalign}
where $M_p(T-\tau)>0$ is independent of $u$ and $t$.\\
\noindent
Now, let $0< \tau<T$. We will bootstrap the $L^1$ bounds in (\ref{eq:2.3}) to better estimates by applying duality arguments. To this end, let $1<q<\infty$ and $\theta \in L^q(\Omega \times (\tau,T))$ such that $\theta \ge0$ and $\|\theta\|_{q,\Omega\times(\tau,T)}=1$. Let $\varphi$  be the unique nonnegative solution in $W_q^{(2,1)}(\Omega\times(\tau,T))$ solving the following system on the left, which at a first glance, appears to be a backward heat equation. But the substitutions $\psi(x,t)=\varphi(x,T+\tau-t)$ and $\tilde\theta(x,t)=\theta(x,T+\tau-t)$ lead to the standard initial boundary value problem given in the right 
\begin{flalign} \label{eq:2.4}
\begin{aligned}[c]
 \begin{cases}
    \varphi_t & = - d_2 \Delta \varphi-\theta, \ \ \ \ \ \ \Omega \times (\tau,T),\\
    \frac{\partial \varphi}{\partial \eta} & = 0, \hspace{0.77 in} \partial\Omega \times (\tau,T),\\
    \varphi & = 0, \hspace{0.85 in} \Omega \times \{T\}.
\end{cases}   
\end{aligned} \ \ \ \ \ \ \ 
\begin{aligned}[c]
\begin{cases}
    \psi_t & = d_2 \Delta \psi+\tilde\theta(x,T+\tau-t), \ \ \ \ \ \ \Omega \times (\tau,T),\\
    \frac{\partial \psi}{\partial \eta} & = 0, \hspace{1.48 in} \partial\Omega \times (\tau,T),\\
    \psi & = 0, \hspace{1.56 in} \Omega \times \{\tau\}.
\end{cases}   
\end{aligned}
\end{flalign}
Results in \cite{ladyzhenskaia1968linear} imply there exists $C_{q,T-\tau}>0$ so that $\|\varphi\|_{q,\Omega\times(\tau,T)}^{(2,1)}\le C_{q,T-\tau}$ and Lemma 4.1 in \cite{morgan1989global} gives $\|\varphi(\tau)\|_{q,\Omega}\le C_{q,T-\tau}$. Now, we start with
\begin{flalign} 
 & \int_{\tau}^T \int_{\Omega} \big(u(x,t)+v(x,t) \big) \theta(x,t) dxdt = \int_{\tau}^T \int_{\Omega} \big(u(x,t)+v(x,t)\big) \big(-d_2\Delta \varphi (x,t)-\varphi_t(x,t)\big) dxdt \nonumber\\
  & = -d_2 \int_{\tau}^T \int_{\Omega}  \big(u(x,t) \Delta \varphi(x,t) + v(x,t) \Delta \varphi(x,t) \big) dxdt - \int_{\tau}^T \int_{\Omega} \big(u(x,t)\varphi_t(x,t) + v(x,t) \varphi_t(x,t) \big) dxdt \nonumber\\
& = A+B \label{eq:2.5}
\end{flalign}
where,
\begin{flalign*}
A & = -d_2 \int_{\tau}^T \int_{\Omega}  (u(x,t) \Delta \varphi(x,t) + v(x,t) \Delta \varphi(x,t)) dxdt \\
& = -d_2 \int_{\tau}^T \int_{\Omega}  u(x,t) \Delta \varphi(x,t) \ dxdt -d_2 \int_{\tau}^T \int_{\Omega} v(x,t) \Delta \varphi(x,t) dxdt \\
& \le d_2 \Vert u(\cdot,t) \Vert_{p,\Omega \times (\tau,T)} C_{q,T-\tau} -d_2 \int_{\tau}^T \int_{\Omega} v(x,t) \Delta \varphi(x,t) dxdt \\
\therefore \ A & \le d_2 M_p(T-\tau) C_{q,T-\tau} -d_2 \int_{\tau}^T \int_{\Omega} v(x,t) \Delta \varphi(x,t) dxdt 
\end{flalign*}
where we have used Holder's inequality and then the estimate from (\ref{eq:2.2}) for $(\Omega \times (\tau,T))$. Now,
\begin{flalign*}
B & = - \int_{\tau}^T \int_{\Omega}  (u(x,t) \varphi_t(x,t) + v(x,t) \varphi_t(x,t)) dxdt \nonumber\\
 & = - \int_{\Omega} \big[\big(u(x,t)+v(x,t) \big)\varphi(x,t) \big]_{\tau}^T dx + \int_{\tau}^T \int_{\Omega} \big(u_t(x,t)+v_t(x,t) \big) \varphi(x,t) dxdt \nonumber\\
 & \le - \int_{\Omega} \big[\big(u(x,T)+v(x,T)\big)\varphi(x,T)-\big(u(x,\tau)+v(x,\tau)\big)\varphi(x,\tau) \big] dx + \int_{\tau}^T \int_{\Omega} \big[\big(d_1 \Gamma u(x,t) + d_2 \Delta v(x,t)\big) \\
 & + f - f \tilde{\kappa}(u(x,t)+v(x,t))\big] \varphi(x,t) dxdt \\
 & \le \int_{\Omega} \big(u(x,\tau)+v(x,\tau)\big)\varphi(x,\tau) dx + \int_{\tau}^T \int_{\Omega} \big(d_1 \Gamma u(x,t)  + d_2 \Delta v(x,t)\big) \varphi(x,t) dxdt + f \int_{\tau}^T \int_{\Omega} \varphi(x,t) dxdt \\
 & - f \tilde{\kappa} \int_{\tau}^T \int_{\Omega} \big(u(x,t)+v(x,t)\big) \varphi(x,t) dxdt 
\end{flalign*}
If $q=\frac{n+4}{2}$ then $W^{(2,1)}_q(\Omega \times (\tau, T))$ embeds continuously into $C(\overline\Omega \times [\tau, T])$ \cite{ladyzhenskaia1968linear}. Using this along with (\ref{eq:2.3}) we get,
\begin{flalign*}
 \int_{\Omega} \big(u(x,\tau)+v(x,\tau)\big) \varphi(x,\tau) dx & \le K_{q,T-\tau} \int_{\Omega} \big(u(x,\tau) + v(x,\tau)\big) dx \le K_{q,T-\tau}M, \\
 \int_{\tau}^T \int_{\Omega} \big(u(x,t)+v(x,t)\big) \varphi(x,t) dxdt & \le K_{q,T-\tau} M(T-\tau) \ \ \text{and } \ \int_{\tau}^T \int_{\Omega} \varphi(x,t) dxdt \le M_{q,T-\tau}
\end{flalign*}
where for the last integral, we have used the Trace Embedding Theorem in \cite{ladyzhenskaia1968linear} to get that there exists $M_{q,T-\tau}>0$ independent of $\theta$ such that $\Vert \varphi \Vert_{L^1,\Omega \times (\tau,T)} \le M_{q,T-\tau}$. Also, (\ref{eq:1.3}) and (\ref{eq:2.2}) give
\begin{flalign*}
    \int_{\tau}^T \int_{\Omega} d_1 \Gamma u(x,t) \varphi(x,t) \ dxdt & = \int_{\tau}^T \int_{\Omega} \int_{\Omega} d_1 \phi(x,y) (u(y,t)-u(x,t)) \varphi(x,t) \ dxdydt \\
    & \le \widehat{K}_{q,T-\tau}
\end{flalign*}
with the continuous embedding of $W^{(2,1)}_q(\Omega \times (\tau, T))$ into $C(\overline\Omega \times [\tau, T])$. Lastly, 
\begin{flalign*}
    \int_{\tau}^T \int_{\Omega} d_2 \Delta v(x,t) \varphi(x,t) \ dxdt  = d_2 \int_{\tau}^T \int_{\Omega} v(x,t) \Delta \varphi(x,t) \ dxdt. 
\end{flalign*}
Replacing these estimates in $B$ and then replacing $A,B$ in  $(\ref{eq:2.5})$ we get,
\begin{flalign*} 
\int_{\tau}^T \int_{\Omega} \big(u(x,t)+v(x,t)\big) \theta(x,t) \ dxdt & \le d_2 M_p(T-\tau) \Vert u(\cdot,\tau) \Vert_{p, \Omega} C_{q,T-\tau} + K_{q,T-\tau} M + \widehat{K}_{q,T-\tau} + f M_{q,T-\tau} \\
& + f \tilde{\kappa} K_{q,T-\tau} M(T-\tau)\\
\therefore \Vert u(\cdot,t) + v(\cdot,t) \Vert_{q',\Omega \times (\tau,T)} & \le N_{q,T-\tau}(1+M)
\end{flalign*}
where, $N_{q,T-\tau} = \max\big\{d_2 M_p(T-\tau) \Vert u(\cdot,\tau) \Vert_{p, \Omega} C_{q,T-\tau} + \widehat{K}_{q,T-\tau} + f M_{q,T-\tau} + f \tilde{\kappa} K_{q,T-\tau} M(T-\tau), K_{q,T-\tau} \big\}$.\\
\noindent
Then by choosing $\tau$ such that $T-\tau$ is sufficiently small and using $q'=\frac{n+4}{n+2}$, we find there exists,
\begin{flalign*} 
\Vert v(\cdot,t) \Vert_{\frac{n+4}{n+2},\Omega \times (\tau,T)} \le N_{q,T-\tau} (1+M)
\end{flalign*}
Note that since this estimate only depends on the size of $T-\tau$, it follows from a limiting process that
\begin{flalign*} 
\Vert v(\cdot,t) \Vert_{\frac{n+4}{n+2},\Omega \times (\tau,T_m)} \le N_{q,T_m-\tau} (1+M)
\end{flalign*}
Now, let us write the above estimate as
\begin{flalign} \label{eq:2.6}
\Vert v(\cdot,t) \Vert_{\frac{n+4}{n+2},\Omega \times (\tau,T_m)} \le K_1
\end{flalign}
Moreover, by Lemma 4.1 from \cite{morgan1989global}, we can obtain the corresponding time–trace estimate
\begin{flalign} \label{eq:time-trace}
    \|v(\tau)\|_{\frac{n+4}{n+2},(\Omega)} \le K_1,
\end{flalign}
\noindent
We claim that if $r \in \mathbb{N}$ then there exists $K_r>0$ such that 
\begin{equation}\label{eq:2.7}
\|v(\cdot,t)\|_{\left(\frac{n+4}{n+2}\right)^r,\Omega\times(0,T_m)}\le K_r
\end{equation}
To this end, note that if $1<q<\infty$ such that $q>\frac{n+2}{2}$ then $W_q^{(2,1)}(\Omega\times(\tau,T))$ embeds continuously into $C(\overline\Omega\times[\tau,T])$ and if $q<\frac{n+2}{2}$ then $W_q^{(2,1)}(\Omega\times(\tau,T))$ embeds continuously into $L^{\alpha}(\Omega\times(\tau,T))$ for all $1\le \alpha \le \frac{(n+2)q}{n+2-2q}$ \cite{ladyzhenskaia1968linear}.  Now, suppose $r\in\mathbb{N}$ such that (\ref{eq:2.7}) is true.\\
\noindent
Now, we want to obtain an estimate as in (\ref{eq:2.6}) with $(\frac{n+4}{n+2})^r$ replaced by $(\frac{n+4}{n+2})^{r+1}$. Let $Q=\left(\frac{n+4}{n+2}\right)^r$, $R=\left(\frac{n+4}{n+2}\right)^{r+1}$ and $q=\frac{R}{R-1}$. If $q\ge\frac{n+2}{2}$ and $\varphi$ solves the left system in (\ref{eq:2.4}), then from above $\|\varphi\|_{p,\Omega\times(\tau,T)}$ is bounded independent of $\theta$ for every $1\le p<\infty$.  On the other hand, if $q<\frac{n+2}{2}$ then $\frac{n+4}{n+2} \le \frac{n+2+2R}{n+2}$ implies that $1 \le \frac{Q}{Q-1}\le\frac{(n+2)q}{n+2-2q}$. As a result, in either case, $\|\varphi\|_{\frac{Q}{Q-1},\Omega\times(\tau,T)}$ is bounded independent of $\theta$. So, calculating again $A,B$ and replacing them in (\ref{eq:2.5}), we can obtain bound as in (\ref{eq:2.7}) with $r$ replaced by $r+1$. Therefore, (\ref{eq:2.7}) is true. Similar argument is also applicable with (\ref{eq:time-trace}).\\
Now, let us recall that
\begin{flalign*}
    g_2(x,t) = u(x,t)v(x,t)^2-(f+\kappa)v(x,t) \le u(x,t)v(x,t)^2
\end{flalign*}
Now, if we consider $G(x,t) = u(x,t)v(x,t)^2$, then $g_2(x,t) \le G(x,t)$. Then (\ref{eq:2.7}) implies $\Vert G(\cdot,t) \Vert_{p,\Omega \times (\tau,T_m)}$ is bounded for each $1<p<\infty$. Now, let $\Psi$ be the unique nonnegative solution to 
\begin{flalign*}
    \Psi_t(x,t) & = d_2 \Delta \Psi(x,t)+G(x,t), \ \ \ \ \ \ \Omega \times (0,T_m),\\
    \frac{\partial \Psi}{\partial \eta} & = 0, \hspace{1.4 in} \partial\Omega \times (0,T_m),\\
    \Psi & = v_0(x) \ge 0, \hspace{0.9 in} \Omega \times \{0\},
\end{flalign*}
then comparison principle for parabolic initial boundary value problems implies $v\le \Psi$, and using the $L^p(\Omega\times(0,T_m))$ bound on $G(x,t)$ for $p>\frac{n+2}{2}$, we observe that $\|\Psi\|_{\infty,\Omega\times(0,T_m)}$ is bounded. $$\Vert v(\cdot,t) \Vert_{\infty,\Omega \times (\tau,T)} \le \tilde{C}$$
which depends on $d_1,d_2,f,\kappa$ and on the domain $\Omega$ but independent of time. This along with (\ref{eq:2.2}) proves the theorem.\hfill $\square$

\section{Numerical Pattern Formation}\label{sec:patterns}
\begin{figure}[H]
\centering
     \includegraphics[width=7cm]{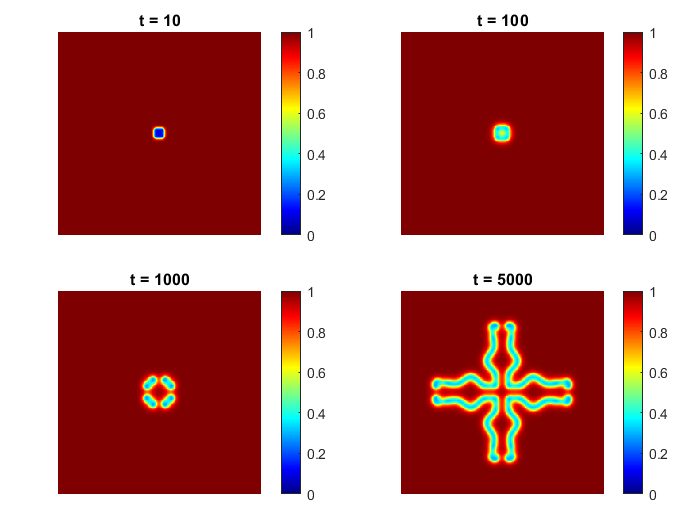}
     \includegraphics[width=7cm]{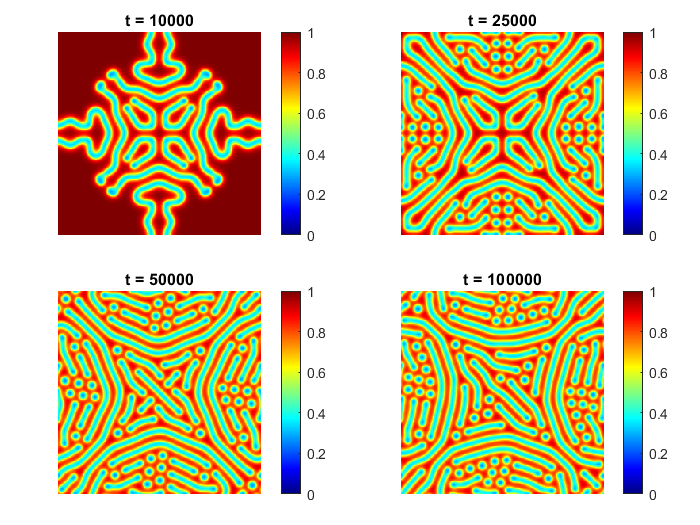}
     \caption{Gray-Scott patterns with local diffusion.}
    \label{fig1: Gray-Scott patterns with local diffusion.}
\end{figure}

\begin{figure}[H]
\centering
     \includegraphics[width=7cm]{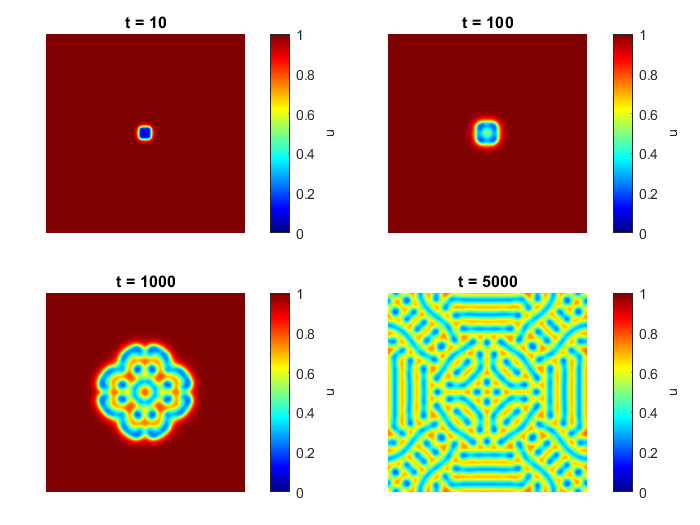}
     \includegraphics[width=7cm]{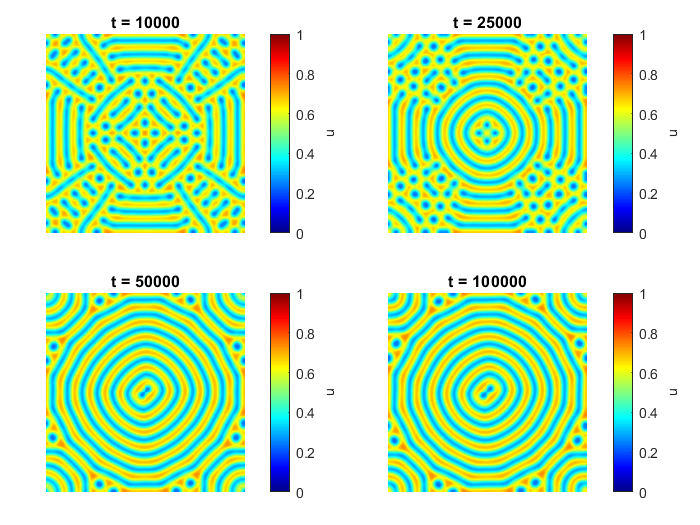}
     \caption{Gray-Scott patterns with local and nonlocal diffusion.}
    \label{fig2: Gray-Scott patterns with local and nonlocal diffusion.}
\end{figure}

The study of pattern formation in reaction--diffusion systems originates with the seminal work of Turing \cite{turing1952chemical}, who showed that diffusion, counterintuitively, can destabilize a spatially homogeneous equilibrium and give rise to spatially heterogeneous steady states. This mechanism, now known as \emph{diffusion-driven instability}, provides a fundamental explanation for the emergence of biological and chemical patterning. Activator--inhibitor systems introduced later by Gierer and Meinhardt \cite{gierer1972theory} offered a concrete biochemical framework in which a slowly diffusing activator promotes its own production while a more rapidly diffusing inhibitor suppresses it, thereby creating spatial contrast.\\
\noindent
Reaction--diffusion systems exhibiting this mechanism are now known to generate a wide variety of organized structures through the interplay between nonlinear reaction kinetics and differential diffusion \cite{murray2003mathematical,ali2023spatiotemporal,pearson1993complex,lee1994pattern,mazin1996pattern,schenk1998interacting}. The Gray--Scott model is a canonical example of such a system, in which the activator $u$ undergoes autocatalytic amplification while simultaneously consuming the inhibitor $v$. The balance between localized growth and diffusion-mediated spreading gives rise to a remarkably rich repertoire of stationary and dynamic patterns, including spot arrays, labyrinths, target waves, spiral interactions, and spatiotemporal chaos; see, for example, \cite{nishiura2001spatio,mcgough2004pattern,muratov2001spike,doelman1998stability}. \\
\noindent
These dynamics are highly sensitive to parameter values and diffusion rates, and small changes can lead to transitions between qualitatively distinct geometric structures. The Gray--Scott model therefore serves as a fundamental testbed for exploring nonlinear pattern formation, both analytically and computationally.\\
\noindent
For the local Gray--Scott system \eqref{eq:local}, the parameters $f$ (feed rate) and $\kappa$ (removal rate) play a central role in pattern selection. In the numerical experiments below, we take
$f = 0.04$, $\kappa = 0.0636$, $d_1 = 1.0$, $d_2 = 0.5$,
a well-studied parameter set known to generate slowly evolving spot and labyrinthine patterns. In this regime, localized concentrations of the inhibitor $v$ nucleate and interact, eventually organizing into persistent spatial structures. We use Central Difference Formula to calculate the local Laplacian. For the system in (\ref{eq:1.2}), we use two dimensional Gaussian Kernel with variance, $\sigma=1$. 
\begin{flalign*}
    \phi(x,y) = \frac{1}{2\pi\sigma^2} \exp\left( -\frac{x^2 + y^2}{2\sigma^2} \right)
\end{flalign*}
All simulations are performed on the periodic square domain $\Omega=[0,L]^2$ with $L = 200$ using Finite Difference Method along with Forward Euler Method. The domain is discretized using a uniform $N \times N$ grid with $N=200$, giving a spatial step size $h = L/N = 1$. The initial conditions are taken to be spatially homogeneous, $u(x,0)=1$, $v(x,0)=0$
with naturally occurring floating-point perturbations serving as symmetry-breaking seeds that initiate pattern formation. We obtain following results. The time step is selected to satisfy the stability requirements of the reaction--diffusion system. The resulting spatiotemporal patterns reveal significant differences between the local and mixed local--nonlocal Gray--Scott models. In particular, Figure~\ref{fig1: Gray-Scott patterns with local diffusion.} -~\ref{fig2: Gray-Scott patterns with local and nonlocal diffusion.} illustrate the qualitative and quantitative distinctions in the pattern dynamics generated by the two diffusion mechanisms.

\section{Conclusion}\label{Conclusion}
In this work, we established the global existence and uniform boundedness of component-wise nonnegative solutions to a Gray--Scott reaction--diffusion system in which the activator diffuses nonlocally while the inhibitor undergoes classical (local) diffusion. The presence of the nonlocal convolution operator introduces a broader and biologically relevant diffusion mechanism that can represent long-range transport, signaling, or dispersal phenomena. Despite this added complexity, the system remains analytically tractable through a combination of semigroup methods and duality techniques, which allow us to derive uniform a priori bounds and prevent finite-time blow-up.\\
\noindent
Our results provide a rigorous foundation for the study of pattern formation in systems that couple local and nonlocal transport processes. Numerical experiments indicate that the introduction of nonlocal diffusion can alter both the geometry and stability of emergent patterns, suggesting that nonlocality plays a significant role in the selection and persistence of spatial structures.\\ 
\noindent
Future research directions include a detailed stability analysis of spatially homogeneous and patterned steady states, characterization of diffusion-driven (Turing) instabilities in the mixed-diffusion setting, and a systematic classification of pattern types as the nonlocal interaction scale varies. Such investigations may contribute to a deeper understanding of morphogenesis, chemical reaction processes, and pattern-forming mechanisms in systems with multiscale transport behavior.\\
\\
\noindent
\textbf{Conflict of Interest:} There is no conflict of Interest.\\
\\
\noindent
\textbf{Acknowledgment:}
I am grateful to my advisor, Dr.\ Jeff Morgan, for his guidance and support throughout this work.

\bibliographystyle{acm}
\bibliography{reff}

\end{document}